\documentclass[11pt]{article}

\usepackage{amsmath,amssymb,amsfonts}

\usepackage{curves}
\usepackage{xypic}  
\usepackage[mathscr]{eucal}
\usepackage[pdftex]{hyperref}

\font\tenmsb=msbm10
\font\sevenmsb=msbm7
\font\fivemsb=msbm5

\newfam\msbfam
\textfont\msbfam=\tenmsb
\scriptfont\msbfam=\sevenmsb
\scriptscriptfont\msbfam=\fivemsb

\font\teneufm=eufm10
\font\seveneufm=eufm7
\font\fiveeufm=eufm5
\newfam\eufmfam
\textfont\eufmfam=\teneufm
\scriptfont\eufmfam=\seveneufm
\scriptscriptfont\eufmfam=\fiveeufm
\def\frak#1{{\fam\eufmfam\relax#1}}

\usepackage[all]{xy}

\newcommand\qed{{\hspace*{\fill}Q.E.D.\vskip12pt plus 1pt}}
\newcommand\Pic[1]{\hbox{\rm Pic(}#1\hbox{\rm )}}

\newcommand\sA{{\cal A}}

\newcommand\sE{{\cal E}}

\newcommand\sL{{\cal L}}
\newcommand\sM{{\cal M}}

\newcommand\sO{{\cal O}}

\newcommand\sQ{{\cal Q}}

\newcommand\sV{{\cal V}}

\newcommand\scr{{\mathscr R}}

\newcommand\scd{{\mathscr D}}

\newcommand\scv{{\mathscr V}}
\newcommand\scy{{\mathscr Y}}

\newcommand\gra{\alpha}
\newcommand\grb{\beta}

\newcommand\vphi{\varphi}

\newcommand\grs{\sigma}

\newcommand\rat{{\mathbb Q}}
\newcommand\reals{{\mathbb R}}
\newcommand\comp{{\mathbb C}}

\newcommand\pn[1]{{\mathbb P}^{#1}}

\newcommand\ft{{\frak t}}

\newcommand\proof{\noindent{\em Proof.}\ \ }

\newtheorem{theorem}{Theorem}[section]
\newtheorem{thm-def}[theorem]{Theorem-Definition}
\newtheorem{lemma}[theorem]{Lemma}
\newtheorem{corollary}[theorem]{Corollary}

\newtheorem{prop}[theorem]{Proposition}

\newtheorem{question}[theorem]{Question}

\newtheorem{definition}[theorem]{Definition}
\newtheorem{re}[theorem]{Remark}
\newtheorem{defre}[theorem]{Definition--Remark}
\newtheorem{pargrph}[theorem]{}
\newtheorem{examp}[theorem]{Example}

\newtheorem{MM}[theorem]{ }

\makeatletter
\ifnum\@ptsize=0\fi
\ifnum\@ptsize=2\fi

\newenvironment{rem*}{\begin{re}\em}{\end{re}}
\newenvironment{example*}{\begin{examp}\em}{\end{examp}}
\newenvironment{definition*}{\begin{definition}\em}{\end{definition}}
\newenvironment{question*}{\begin{question}\em}{\end{question}}

\newenvironment{prgrph*}[1]{\indent\begin{pargrph}{\bf #1.}\em\
}{\end{pargrph}}
\newenvironment{defre*}{\begin{defre}\em}{\end{defre}}
\newenvironment{MM*}{\begin{MM}\em}{\end{MM}}

\begin{document}

\title{Geometry  of rays-positive manifolds
}
\author{M.C. Beltrametti, A.L. Knutsen, A. Lanteri, and C. Novelli}

\date{}

\maketitle


\begin{abstract} Let $\sM$ be a smooth complex projective variety  and let $\sL$ be a line bundle on it. 
Rays-positive manifolds, namely pairs $(\sM,\sL)$ such that $\sL$ is numerically effective and $\sL\cdot R>0$ for all extremal rays $R$ on $\sM$, are studied. Several illustrative examples and some applications are provided. In particular, projective  varieties with   crepant singularities and  of small degree with respect to the codimension are classified, and the non-negativity of   the sectional genus $g(\sM,\sL)$ is proven, describing as well the  pairs with $g(\sM,\sL)=0,1$.

\end{abstract}

\section*{Introduction} Let $\sM$ be a  smooth complex projective variety of dimension $n\geq 2$, and let $\sL$ be a line bundle on $\sM$. 
Assume that $K_\sM$ is not numerically effective (nef).
In classical adjunction theory    $\sL$ is assumed to be ample.  Then, by the Kawamata rationality theorem, the invariant
$$\tau=\tau(\sM,\sL) := \inf\{t \in \reals\ |\ K_\sM + t\sL \ \mbox{is  nef} \}$$ is a positive rational number, the {\em nefvalue} of $(\sM,\sL)$.  The classical adjunction theoretic approach to   the classification of  polarized manifolds $(\sM,\sL)$  is based on the study of  the structure of the  morphism  associated to  the divisor $K_\sM+\tau \sL$ (see \cite{Book}).

One main obstruction to extending this study to the case when $\sL$ is merely nef   
 is given by the possible existence of  cycles $Z \in \overline{NE}(\sM)$ such that $K_\sM\cdot Z<0$ and $\sL\cdot Z=0$. In this case, the invariant $\tau$ is not defined. To overcome this problem,  for any extremal 
ray $R=\reals_+[C]$, with $C$ a minimal rational curve, such that $\sL\cdot C>0$ (such an extremal ray will be called $\sL$-{\em positive}), we define the invariant
$$\tau_\sL(R):= \frac{-K_\sM\cdot C}{\sL\cdot C}$$
(see Definition \ref{CA}). This  does not require $\sL$ to be nef, so we can in fact work with any  line bundle $\sL$, that is, with  any {\em pre-polarized manifold}  $(\sM,\sL)$. Let  $\vphi\colon\sM \to Y$ be  the 
  contraction associated to an $\sL$-positive extremal ray $R$. 
Since $\sL$ is $\vphi$-ample, there exists  an ample line bundle $A$ on $Y$ such that $\sL + \vphi^*A$ is ample on $\sM$. Clearly, $\tau_\sL(R)=\tau_{\sL+\vphi^*A}(R)$, and   the invariant $\tau_\sL(R)$ is just the nefvalue of the polarized variety $(\sM, \sL+\vphi^*A)$. 

In  Section \ref{SF}  we recall  some
 structure results we need about  pairs $(\sM,\sL)$ admitting  an $\sL$-positive extremal ray $R$. As noted above,  they follow from  the corresponding classification results of extremal rays and Fano--Mori contractions in the case of {\it polarized} manifolds. An iterative application of these results 
leads as well  to a  natural definition of  a first reduction map  for any pre-polarized manifold $(\sM,\sL)$ in terms of local contractions.

  In Section \ref{AMP}  we introduce the notion of rays-positive manifold $(\sM,\sL)$.
  If    $\sL$ is nef,  one can define the invariant
$\frak{t}=\frak{t}(\sM,\sL) := \sup\{t \in \reals\ |\ tK_\sM + \sL \ \mbox{is  nef} \}$,  which is a non-negative rational number, again by the Kawamata rationality theorem.  We define 
 $\sL$ to be {\em rays-positive}, and $(\sM,\sL)$ to be a {\em rays-positive manifold}, if $\frak{t}(\sM,\sL) >0$ (so that  the number $\frac{1}{\frak{t}(\sM,\sL)}$     corresponds  to the nefvalue of $(\sM,\sL)$ in the polarized case). Our terminology comes from the fact  that  $\sL$ is rays-positive if and only if  all extremal rays are 
$\sL$-positive (see Lemma \ref{Dm}). The class of rays-positive manifolds $(\sM,\sL)$  is  somehow  the largest class of pre-polarized manifolds  where adjunction still works, in the sense that $K_\sM+k\sL$ will be nef (whence semi-ample if $\sL$ is big) for large $k$. 

Recalling that a line bundle $\sL$ is  said to be  {\em numerically positive} ({\em nup}) if $\sL\cdot C> 0$ for all  curves $C$ on $\sM$,   one  has that nup line bundles are rays-positive. (Nup line bundles are also called {\em strictly nef} in the literature, see e.g., \cite{Serrano}.)  We  show that the converse is not true (see Examples  \ref{Criscito}, \ref{quasiF}, \ref{AnCa}).  This fact together with  several existing  results  on nup line bundles (see e.g., \cite{Serrano}, \cite{IoNef}) give a further  motivation to study rays-positive manifolds.

Another   motivation comes from the case of surfaces. It is easily seen  that a {\em quasi-polarized } surface $(\sM,\sL)$ (that is, $\sL$ is nef and big) is rays-positive if and only if there are no $(-1)$-curves $E$ on $\sM$ such that  $\sL\cdot E=0$. According to standard terminology, such a pair  $(\sM, \sL)$ is called $a$-{\em minimal} (or $\sL$-{\em minimal}). 
There exists a wide literature on such surfaces: we refer in particular to \cite[\S 7]{CR}, which also  includes results for  $\sL$  merely nef. 

A  relevant geometric  context  where rays-positive pairs occur is the case  of a projective variety $X$ with crepant singularities:  indeed  if $\pi:\widetilde{X}\to X$ is any crepant resolution, then $(\widetilde{X}, \pi^*\sO_X(1))$ is rays-positive (cf. Definition \ref{def:crepsing} and Lemma \ref{crep-am}).

In  Sections \ref{FRPM}, \ref{lowdeg} and \ref{FC} we provide some further results and applications for rays-positive manifolds. Let us briefly mention some of them.

 Rays-positive quasi-polarized manifolds $(\sM,\sL)$ such that $K_\sM+(n-1)\sL$ is not nef and big are classified up to first reduction in Corollary \ref{9ago}. In particular, they are all uniruled of $\sL$-degree at most one with the exception of $(\pn 3, \sO_{\pn 3}(2))$, cf. Corollary \ref{cor:nonid-rp}.  Similarly,  rays-positive  manifolds $(\sM,\sL)$ such that $K_\sM+(n-2)\sL$ is not pseudo-effective  are classified in Proposition \ref{prop:class}. The latter result leads  in Theorem \ref{thM:lowdeg} to a classification  of   projective  varieties with   crepant singularities and  of small degree with respect to the codimension. 
This generalizes a result obtained by Ionescu \cite{IoSmall} in the setting of smooth varieties. In Example \ref{Crecon} we construct  such varieties with crepant singularities.
 
 In propositions \ref{genus} and \ref{g01}, we show that $g(\sM,\sL)\geq 0$ for a rays-positive  manifold $(\sM,\sL)$ of any dimension, regardless the bigness of $\sL$ (this was conjectured by Fujita \cite{FujitaQ}   for any quasi-polarized normal variety $(\sM,\sL)$, and proved  in dimension $\leq 3$) 
and we  describe rays-positive manifolds of sectional genus $g(\sM,\sL)=0,1$.
The scroll of sectional genus one described in Example  \ref{Ascroll} also shows that the
 inequality $g(\sM,\sL)\geq h^1(\sO_\sM)$, conjectured in the setting of quasi-polarized varieties, is not true dropping the bigness assumption.

\smallskip

{\bf Notation and terminology.} We work on the complex field $\comp$ and use the standard terminology in
algebraic geometry. In particular, 
we use the additive notation for  the tensor product of line bundles on a projective variety $X$, and by $K_X$ we denote the canonical bundle if  $X$ is smooth.

If  $X$ is smooth and $L$ is  {\em any} line bundle on  $X$, we say that the pair $(X,L)$  is a {\em pre-polarized manifold}.
The {\em sectional genus} $g(X,L)$ of  $(X,L)$ is defined by 
$2g(X,L)-2= (K_X+(n-1)L)\cdot L^{n-1}$. It is well-known that  $g(X,L)$ is an integer, cf. e.g., \cite[p. 25]{FuBook}.

A pre-polarized manifold $(X,L)$ is called a {\em
scroll} over a smooth $m$-dimensional variety $Y$ if there is a surjective morphism $\pi: X\to Y$ such
that $(F,L_F)\cong (\pn {n-m},\sO_{\pn {n-m}}(1))$ for every fiber $F$.  (We allow the case $m=0$.)
Since $L$ is $\pi$-ample, we have that $L+\pi^*A =:\sA$ is ample for some  very ample line bundle $A$ on $Y$ (see  \cite[Proposition 1.45]{KM}) and $\sA_F\cong\sO_{\pn {n-m}}(1)$ for each fiber $F$ of $\pi$. Therefore $(X,\sA)\cong ({\mathbb P}(\sV),\xi_\sV)$, where $\xi_\sV$ is the tautological line bundle on $X$ of the ample vector bundle $\sV:=\pi_*\sA$
of rank $n-m+1$ on $Y$ (see e.g.,  \cite[Proposition 3.2.1]{Book}). We then have  $\sE:=\sV\otimes (-A)=\pi_*L$, so that $(X,L)\cong ({\mathbb P}(\sE),\xi_\sE)$.
Also note that $g(X,L)=g(Y)$ if $Y$ is a curve, by the Chern--Wu relation.

According to \cite{FujitaQ}, a pre-polarized manifold $(X,L)$ is said to be a {\em quasi-polarized manifold} if  $L$ is  nef and big.
We say that $X$ is a {\em quasi-Fano} manifold if  $ -K_X$ is nef  and big. 
Notice that quasi-Fano manifolds are often called {\em almost Fano manifolds}, as well as {\em weak Fano manifolds} in the literature. 

A quasi-polarized manifold $(X,L)$ is called a {\em quasi-Del Pezzo manifold} (resp., a {\em quasi-Mukai manifold}) if $-K_X=(n-1)L$ (resp., $-K_X=(n-2)L$ with $n\geq 3$). If $L$ is ample, the prefix ``quasi'' is deleted. (Note that, according to our terminology, Del Pezzo manifolds with $\varrho\geq 2$, where $\varrho$ denotes the Picard number,  are also scrolls over surfaces for $n\geq 3$, cf. \cite{FuBook}).

\section{$\sL$-positive extremal rays}\label{SF}
\addtocounter{subsection}{1}\setcounter{theorem}{0}

 In this section we collect results that we will need in the rest of
the paper. 

Let $\sM$ be a smooth projective variety and let $\sL$ be a line bundle on $\sM$. If $K_\sM$ is not nef, it is well known that there exists (at least) an extremal ray on $\sM$. We will always write an extremal ray $R$ as $R=\reals_+[C] $, where $C$ is a rational curve of minimal anticanonical degree among curves whose numerical class belongs to $R$, and we will denote the  length of $R$ by $\ell(R) := -K_\sM \cdot C$.

\begin{definition*}\label{CA} Let $\sM$ be a smooth projective variety and let $\sL$ be a line bundle on $\sM$. We say that an extremal ray 
$R=\reals_+[C] $ on $\sM$ is $\sL$-{\em positive} if $\sL\cdot C>0$. For such a ray  set (cf. \cite{IoNef})
$$\tau_\sL(R):=\frac{\ell(R)}{\sL\cdot C}.$$
\end{definition*}

Note that an extremal ray $R$ is orthogonal to a given  adjoint bundle $tK_\sM+\sL$, where $t$  is a positive constant, if and only if  $ \tau_\sL(R)=1/t$.

Now let $(\sM,\sL)$ be a pre-polarized manifold of dimension $n \geq 2$ and $\vphi\colon\sM \to Y$  the 
  contraction associated to an $\sL$-positive extremal ray $R$. 
Since $\sL$ is $\vphi$-ample, there exists  an ample line bundle $A$ on $Y$ such that $\sL + \vphi^*A$ is ample on $\sM$ (see  \cite[Proposition 1.45]{KM}). Clearly, $\tau_\sL(R)=\tau_{\sL+\vphi^*A}(R)$; moreover,  $\sL$ and $\sL + \vphi^*A$ are isomorphic on the fibers of $\vphi$. In other words, the invariant $\tau_\sL(R)$ is just the nefvalue of the polarized variety $(\sM, \sL+\vphi^*A)$. 
Therefore classification results of extremal rays and Fano--Mori contractions in the case of {\it polarized} manifolds yield structure results about  pairs $(\sM,\sL)$ admitting  an $\sL$-positive extremal ray $R$ (see \cite{ABW2}, \cite{AO}, \cite{AW}, \cite{AWP}, \cite{Wis}, \cite{WisCont}  and \cite{Book}).

The following  results deal with  all the cases with $\tau_\sL(R)> \dim \sM-2$ we need in the sequel. They are natural extensions of the classical adjunction
theoretic knowledge in the case of ample line bundles \cite{Book}, obtained using contractions of extremal rays instead of the nefvalue morphism.

\begin{prop}\label{highnef} 
 Let $(\sM,\sL)$ be a pre-polarized manifold of dimension $n \geq 2$,  and  let  $R$ be  an $\sL$-positive extremal ray. Then   $\tau_\sL(R)\leq n-1$  unless either
\begin{enumerate}
\item[$(1)$] $\tau_\sL(R) = n + 1$ and $(\sM,\sL)\cong(\pn n ,\sO_{\pn n}(1))$; or
\item[$(2)$] $\tau_\sL(R) = n$ and $(\sM,\sL)\cong (\sQ,\sO_\sQ(1))$,  where $\sQ$ is a smooth hyperquadric in $\pn {n+1}$; or
\item[$(3)$] $\tau_\sL(R) = n$ and $(\sM, \sL)$ is a scroll over a smooth curve $Y$; or
\item[$(4)$]  $n=2$, $\tau_\sL(R)=3/2$ and $(\sM,\sL)\cong(\pn 2,\sO_{\pn 2}(2))$.
\end{enumerate}
\end{prop}

\begin{prop}\label{notbig}  Let $(\sM,\sL)$ be a pre-polarized manifold of dimension $n \geq 2$,  and  let  $R$ be  an $\sL$-positive extremal ray. 
 Assume  $\tau_\sL(R)=n-1$ and let $\vphi:\sM\to Y$  be the contraction associated to $R$. Then one of the following cases  occurs:
\begin{enumerate}
\item[$(1)$]   $(\sM,\sL)$ is a Del Pezzo manifold of Picard number one.
\item[$(2)$] The variety $Y$ is  a  smooth curve, and $(F,\sL_F)\cong(\sQ,\sO_{\sQ}(1))$, with $\sQ$ a reduced and irreducible hyperquadric in $\pn n$,  for every fiber $F$ of $\vphi$, and the general fiber is smooth $($if $n=2$, this means that $(F,\sL_F)\cong(\pn 1,\sO_{\pn 1}(2)))$. 
\item[$(3)$] The variety $Y$ is a  smooth surface and $(\sM,\sL)$ is a scroll over $Y$.
\item[$(4)$] The morphism  $\vphi$ is birational and   contracts to a smooth point  a divisor $E\cong \pn {n-1}$ such that $\sO_E(E)= \sO_{\pn {n-1}}(-1)$, $\sL_E=\sO_{\pn {n-1}}(1)$, and 
$\sL=\vphi^* L -E$, where $L:=(\vphi_*\sL)^{**}$, the double dual. Moreover, $K_\sM+(n-1)\sL=\vphi^*(K_Y+(n-1)L)$.\end{enumerate}

Furthermore, let $\{R_i\}_{i\in I}$ be the family of all non-nef $\sL$-positive extremal rays such that $\tau_\sL(R_i)=n-1$,  and let  $E_i$  be the locus of $R_i$.  If $n\geq 3$, then the exceptional divisors $E_i$ are pairwise disjoint.
\end{prop}

\begin{prop}\label{n21}  Let $(\sM,\sL)$ be a pre-polarized manifold of dimension $n \geq 3$,  and  let  $R$ be  an $\sL$-positive extremal ray. 
 If $n-2<\tau_\sL(R)<n-1$, then either
\begin{enumerate}
\item[$(1)$] $n=4$, $\tau_\sL(R)= 5/2$, and $(\sM,\sL) \cong (\pn 4,\sO_{\pn 4}(2))$; or
\item[$(2)$] $n=3$, $\tau_\sL(R) = 3/2$, and $(\sM,\sL) \cong (\sQ,\sO_\sQ(2))$, $\sQ$ a hyperquadric in $\pn 4$; or
\item[$(3)$] $n=3$, $\tau_\sL(R) = 4/3$, and $(\sM,\sL) \cong (\pn 3,\sO_{\pn 3}(3))$; or
\item[$(4)$] $n=3$, $\tau_\sL(R) = 3/2$,  the  contraction $\vphi\colon \sM\to Y$  associated to $R$  maps onto a smooth curve $Y$, and  $(F,\sL_F) \cong (\pn 2,\sO_{\pn 2}(2))$ for every fiber $F$ of $\vphi$.
\end{enumerate}
\end{prop}

The results above allow us to define a first reduction map for arbitrary pre-polarized manifolds. The definition naturally  extends the classical notion of  first reduction in the adjunction theoretic sense given in the ample case (see e.g., \cite[Chapter 7]{Book}). 

The key observation is that, by  propositions \ref{highnef} and \ref{notbig}, any $\sL$-positive extremal ray $R$ on an $n$-dimensional pre-polarized  manifold $(\sM,\sL)$ with $\tau_\sL(R)\geq n-1$ is nef, except precisely  for the case in Proposition \ref{notbig}$(4)$, where $\tau_\sL(R)=n-1$. In the latter case, the contraction of the ray is birational onto a smooth manifold. In fact, if $\dim \sM \geq 3$, the extremal rays in question are disjoint, so that there is a simultaneous contraction of all such rays 
$\varphi: \sM \to \sM_1$, which is birational and $\sM_1$ is smooth. If $\dim \sM =2$, we can pick a maximal subset of pairwise disjoint such rays (which correspond to $(-1)$-curves) and obtain a similar simultaneous contraction.  We can then repeat the procedure with  the pair $(\sM_1, \sL_1)$, where $\sL_1:= (\varphi_*\sL)^{**}$ is the double dual.
Iterating this process, we obtain at the end a birational morphism $\Phi:\sM \to M$, where $M$ is a smooth projective variety, $\Phi$  is a sequence of contractions of $\pn {n-1}$'s to smooth points, and $M$ does not contain any non-nef $L$-positive extremal ray $\scr$ with $\tau_L(\scr)=n-1$. If $\dim \sM \geq 3$, then  the map $\Phi$ is uniquely determined. In the case $\dim \sM =2$, the map $\Phi$ depends on a choice of which $(-1)$-curves to contract, cf. Example \ref{antonio} below. 

Summarizing, we obtain the following result and definition of  first reduction.

\begin{thm-def}\label{Crespo} {\rm (First reduction for pre-polarized manifolds)}
 Let $(\sM,\sL)$ be a pre-polarized manifold of dimension $n \geq 2$.
 Then there exists  a birational morphism $\Phi:\sM \to M$ onto a smooth projective variety $M$ 
such that $\Phi^*(K_M+ (n-1)L)=K_\sM+(n-1)\sL$, where $L=(\Phi_*\sL)^{**}$ is the double dual, and $M$ does not contain any non-nef $L$-positive extremal ray $\scr$ with $\tau_L(\scr)=n-1$. The morphism $\Phi$ is a sequence of contractions of $\pn {n-1}$'s to smooth points and is uniquely determined, up to isomorphisms, if $n\geq 3$.

We say that  the pair $(M,L)$ and the map $\Phi$ are a {\em first reduction} and  a  {\em first reduction map} of  the pre-polarized manifold $(\sM,\sL)$, respectively.
\end{thm-def}

\smallskip

Unlike the classical case of  polarized manifolds, the first reduction map in the pre-polarized case is not necessarily just a simultaneous contraction of disjoint extremal rays $R$ with $\tau_\sL(R)=n-1$. However, if one requires $\sL$ to be ample or nef, one can easily describe the exceptional locus of $\Phi$. The proof of this fact is an almost straightforward study of local contractions, so we omit it.

\begin{prop}\label{palacio}
 With the same assumptions and notation as in Theorem-Definition \ref{Crespo}, assume   $\sL$ to be  nef $($resp., ample$)$.  Then the exceptional locus of $\Phi$, if non-empty,   consists of disjoint chains of  strict transforms of $\pn {n-1}$'s $($resp.,  disjoint $\pn {n-1}$'s$)$. 
\end{prop}

Recall that classically the notion of reduction is given only for polarized manifolds $(\sM,\sL)$ such that $K_{\sM}+(n-1)\sL$ is nef and big \cite[p.~171]{Book}, and for these pairs our definition coincides with the classical one. However, our definition applies in particular to all polarized manifolds regardless $K_{\sM}+(n-1)\sL$ is nef and big or not. For polarized manifolds, Proposition \ref{nonid-am} below will show that it is in a way  just a trivial extension, in the sense that $\Phi$ is an isomorphism except for a few explicitly described cases.

As  a  consequence of Theorem-Definition \ref{Crespo}, any $L$-positive extremal ray $\scr$ on $M$ with $\tau_L(\scr)\geq n-1$ is necessarily nef and $(M,L)$, as well the contraction of $\scr$, is as in one of propositions \ref{highnef} and \ref{notbig}$(1)$--$(3)$.
In particular, 
let us stress the fact that according to our definition
the reduction $(M,L)$ of a pair $(\sM, \sL)$ might be
covered by lines (that is, smooth rational curves $\ell \subset M$
such that $L \cdot \ell=1$). For instance, see Example \ref{AnCa}. This cannot happen
in the classical case \cite[Theorem 7.6.6$(1)$]{Book}.

We conclude this section  with   an example which shows that the first reduction map is not uniquely determined  when  $n=2$.
\begin{example*} \label{antonio}
  Let $M$ be a $\pn {1}$-bundle over a smooth curve of positive genus and let $L$ be any 
 line bundle satisfying $L \cdot F=2$, where $F$ is the algebraic equivalence class of the fibers. Pick any point $x$ in a fiber $F_0$ and let $\sigma: \sM \to M$ be the blowing-up at $x$. Let $E$ be the exceptional curve and $\widetilde{F}_0$ the strict transform of $F_0$.  Set $\sL:= \sigma^*L-E$. Then one easily sees that
both $E$ and $\widetilde{F}_0$ are $(-1)$-curves satisfying $E \cdot \sL = \widetilde{F}_0 \cdot \sL=1$ and they intersect in one point. In fact, these two curves generate the only two extremal rays on $\sM$. Both extremal rays are as in case $(4)$ of Proposition \ref{notbig}. 
Now we can choose either to contract $E$ or $\widetilde{F}_0$, which leads us to two possible  first reduction maps $\Phi$.

Choosing $\Phi$ to be the contraction of $E$, we have
$\Phi=\sigma: \sM \to M$ and we get back $(M, L)$, which is as in Proposition \ref{notbig}(2). 

Choosing $\Phi: \sM \to M'$ to be the contraction of $\widetilde{F}_0$ we obtain a different pair $(M',(\Phi_*\sL)^{**})$, which, however,  is still as in Proposition \ref{notbig}(2).

Now consider instead the line bundle $\sL':= \sigma^*L$ on $\sM$. The two $(-1)$-curves $E$ and  $\widetilde{F}_0 $ satisfy $E \cdot \sL' = 0$ and $\widetilde{F}_0 \cdot \sL'=2$. Therefore
only the extremal ray $R:=\reals_+[\widetilde{F}_0]$, of length $\ell(R)=1$,  is $\sL'$-positive, with invariant   $\tau_{\sL'}(R)=\frac{1}{2} < n-1=1$. In conclusion, there is 
no $\sL'$-positive extremal ray $R$ on $\sM$ with $\tau_{\sL'}(R)\geq n-1=1$, 
so that  the  first reduction map  with respect to $(\sM,\sL')$ is an isomorphism.\end{example*}

\section{Rays-positive manifolds}
\label{AMP}\addtocounter{subsection}{1}\setcounter{theorem}{0}

In this section we introduce the notion of rays-positive manifold providing first results and several examples.
Let $\sM$ be a smooth projective variety and 
let $\sL$ be a nef line bundle on $\sM$. Assume that $K_\sM$ is not nef and let
$ \ft(\sM,\sL) := \sup\{t \in \reals\ |\ tK_\sM + \sL \mbox{\ is  nef\;}\}$.
Then, by  a  version of the   Kawamata rationality theorem   in the case of  a nef line bundle (see   \cite[Exercise 6.7.5, p. 166]{Debarre} and  \cite[10-3-4]{Mat}),
  $\ft(\sM,\sL)$ is a non-negative rational number. 
Moreover, there is an 
 extremal ray $R$ in $\overline{NE}(\sM)$ such that
$\big(\ft(\sM,\sL)K_\sM+\sL\big)\cdot R=0$.

Note that   $\ft(\sM,\sL)>0$ if $\sL$ is ample,  in which case 
$\ft(\sM,\sL)$ is the reciprocal of the nefvalue of $(\sM,\sL)$ recalled in the introduction.
We now give examples with $\ft(\sM,\sL)=0$. 
\begin{example*}\label{AntoEx0} Consider the $\pn {n-1}$-bundle over $\pn 1$,
$\sM = {\mathbb P}_{\pn 1}(\sV)$, where $\sV =
\bigoplus_{i=1}^n \sO_{\pn 1}(a_i)$ is normalized as in
\cite[Lemma 3.2.4]{Book}, i.e., $a_1 \geq a_2 \geq \dots \geq
a_n=0$. Let $p:\sM \to \pn 1$ be the projection and set
$\sL: = p^* \sO_{\pn 1}(1)$. Clearly $\sL$ is nef but
not big, and $K_{\sM} = -n\xi + p^*\sO_{\pn 1}(a-2)$,
where $\xi$ is the tautological line bundle and $a = \sum_{i=1}^n
a_i$. Then $K_{\sM}$ is not nef and $tK_{\sM}+\sL = -tn\
\xi +  p^*\sO_{\pn 1}\big(t(a-2)+1\big)$ is nef if and only
if $t \leq 0$ and $t(a-2)+1 \geq 0$ \cite[Lemma 3.2.4]{Book}.
Therefore $\frak{t}(\sM, \sL)=0$.
\end{example*}
\begin{example*}\label{Aex0}  Let $\grs\colon \sM\to M$ be  the blowing-up of a smooth $n$-fold $M$ at a point, and let $\sL=\grs^*L$, where
$L$ is an ample line bundle on $M$. For any curve $C$ on $\sM$ contained in  the exceptional divisor $E$ we have $(\sL+tK_\sM)\cdot C \leq 0$ for $t\geq 0$. Thus $\ft(\sM,\sL)=0$.

More generally, we have the following. 
Let $\sM$ be a smooth projective variety, and  let $\vphi\colon \sM\to V$ be a proper birational morphism  where  $V$ is a normal  variety with $\rat$-factorial singularities. Then every irreducible component of  the exceptional locus ${\rm Exc}(\vphi)$ of $\vphi$ 
 has codimension one in $\sM$. Furthermore, there exists an effective $\rat$-Cartier divisor $J$ on $\sM$, whose support is ${\rm Exc}(\vphi)$, and $J\cdot C<0$  for any curve $C$ contracted by $\vphi$ (see e.g., \cite[\S 1.10, p. 28]{Debarre}).  Moreover,  if $V$  has terminal singularities, the equality
$K_\sM=\vphi^*(K_V)+\lambda J$ holds true in $\Pic\sM\otimes\rat$, for some positive rational coefficient  $\lambda$.
Let now $\sL:=\vphi^* H$ for some ample line bundle $H$ on $V$. Then $\sL$ is nef and, for any curve $C$ contracted by $\vphi$, one has
$(tK_\sM +\sL)\cdot C = t  \lambda (J\cdot C)<0$ for each positive $t\in \reals$. This implies  $\ft(\sM,\sL)=0$.
\end{example*}

\begin{definition*}\label{A}  Let $\sM$ be a smooth projective variety and  $\sL$  a nef line bundle on $\sM$. We say that  $(\sM,\sL)$ is a {\em rays-positive manifold}, and that $\sL$ is {\em rays-positive},  if  either $K_{\sM}$ is nef or if $K_{\sM}$ is not nef and
$\ft(\sM,\sL)>0$. In the latter case, $\tau_\sL(R)=\frac{1}{\ft(\sM,\sL)}$ for any extremal ray $R$ orthogonal to $\ft(\sM,\sL) K_\sM+\sL$.
\end{definition*}

\begin{rem*} \label{remnefvalue}  Let $\varphi: \sM \to Y$ be the extremal contraction of the ray $R$ in Definition \ref{A}. As noted in the beginning of \S\ref{SF}, there exists  an ample line bundle $A$ on $Y$ such that $\sL + \vphi^*A$ is ample and $\tau_\sL(R)=\tau_{\sL+\vphi^*A}(R)$. Hence  $\ft(\sM,\sL)= \ft(\sM,\sL+ \vphi^*A)$, the reciprocal of the nefvalue of the
polarized pair $(\sM,\sL+ \vphi^*A)$ in the classical adjunction theoretic sense.
\end{rem*}

By our  definition  all  nef line bundles are  rays-positive whenever  the canonical bundle is nef.  We have included this case for  technical reasons. The relevant framework is clearly when the canonical bundle is non-nef. The terminology  is clarified by the following.

\begin{lemma}\label{Dm}  Let $\sM$ be a smooth projective variety of dimension $n\geq 2$ and let $\sL$ be a nef line bundle on $\sM$.  The following conditions are equivalent:
 \begin{enumerate}
\item[$(1)$] $(\sM,\sL)$ is a rays-positive manifold.
\item[$(2)$] All extremal rays on $\sM$ are $\sL$-positive.
\item[$(3)$] $\sL\cdot C>0$ for all curves $C$ such that $K_\sM\cdot C <0$.
\end{enumerate}

Moreover, if $\sL$ is big, the above are further  equivalent to
\begin{enumerate}
\item[$(4)$]  All non-nef extremal rays on $\sM$ are  $\sL$-positive.
\end{enumerate}
\end{lemma}
\proof  The equivalence between $(1)$ and $(2)$  is an immediate consequence of  the existence of an extremal ray orthogonal to $\frak{t}(\sM,\sL)K_\sM+\sL$, while
the equivalence between $(2)$ and $(3)$  follows from the Mori cone theorem.

Obviously, $(2)$ implies $(4)$. If $\sL$ is big, the converse follows since any
nef extremal ray is $\sL$-positive. Indeed, write $m\sL = A +D$ with $A$ ample and $D$ effective (see \cite[Lemma 2.60(2)]{KM}), and pick a generator $C$ of the ray that is not contained in $D$; then $\sL\cdot C=\frac{1}{m}(A \cdot C + D \cdot C)>0$.
  \qed

The above lemma implies  that a nup line bundle is rays-positive. The converse is not true, as shown by the following  example,  as well as in examples \ref{quasiF} and  \ref{AnCa} below.
(We mention  \cite{Serrano}, \cite{IoNef}, \cite{BSnup} and \cite{CCP} for results on nup line bundles.) 

\begin{example*}\label{Criscito} Let $\sM$ be a $\pn 1$-bundle of positive invariant over a smooth curve of positive genus. Let $E$ be the section with minimal self-intersection $E^2=-e$, with $e> 0$, and let $f$ be a fiber. Take $\sL=a(E+ef)$, $a>0$. Then $\sL$ is nef and big but not nup since $\sL\cdot E=0$. On the other hand $(\sM,\sL)$ is rays-positive according to Lemma \ref{Dm} since the only  extremal ray is $R=\reals_+[f]$. Moreover  one has $(aK_\sM+2\sL)\cdot R=0$, that is, $\ft(\sM,\sL)=\frac{a}{2}$. \end{example*}

Remaining in  the case of surfaces, since all non-nef extremal rays  on a surface are generated by $(-1)$-curves, Lemma \ref{Dm}  says that  a smooth quasi-polarized surface $(\sM,\sL)$ is rays-positive if and only if  there are no $(-1)$-curves $E$ on $\sM$ satisfying $\sL\cdot E=0$, that is, $(\sM,\sL)$ is $a$-{\em minimal}, or $\sL$-{\em minimal}, according to standard terminology in the literature.

We refer to \cite[\S 7]{CR} for an extended study of quasi-polarized surfaces $(\sM,\sL)$, including results where $\sL$ is assumed to be merely nef.

The following  two results show  how  the  invariant $\ft$ and the concept of rays-positivity behave under first reduction.

\begin{lemma}\label{tsale} Let $\sM$ be a smooth projective variety of dimension $n\geq 2$ with $K_\sM$ not nef,  and let $\sL$ be a nef line bundle on $\sM$. Let $(M,L)$ be a first reduction  of $(\sM,\sL)$. Then
$\ft(M,L)\geq \ft(\sM,\sL)$  if $K_M$ is not nef. Moreover, $(\sM,\sL)\cong (M,L)$ if $\ft(\sM,\sL) >\frac{1}{n-1}$.

In particular, if $(\sM,\sL)$ is rays-positive, then  so is $(M,L)$.

\end{lemma}

\proof If $\ft(\sM,\sL)>\frac{1}{n-1}$, then any extremal ray $R$ on $\sM$ satisfies $\tau_\sL(R)\leq \frac{1}{\ft(\sM,\sL)}<n-1$, so that $(\sM,\sL)\cong(M,L)$ by construction of the first reduction map $\Phi$.

If $\ft(\sM,\sL)=0$, there is nothing left to prove, so we can assume $0< \ft(\sM,\sL) \leq\frac{1}{n-1}$.  Now $\ft(M,L)\geq\ft(\sM,\sL)$ follows since $K_M+\frac{1}{\ft(\sM,\sL)}L$ is nef, which is a consequence of the facts that $\Phi^*(K_M+(n-1)L) =K_\sM+(n-1)\sL$ and $\ft(\sM,\sL) \leq\frac{1}{n-1}$.

The last assertion is now  clear.
\qed

The following example involves some of the concepts above. It also illustrates the 
iterative 
procedure behind Theorem-Definition \ref{Crespo}.
\begin{example*}\label{AntoEx1}
 Consider $(\pn 2, \sO_{\pn 2}(2))$. Clearly,
$K_{\pn 2}+\sO_{\pn 2}(2)$ is not nef. Let $\sigma_1: \scy_1 \to
\pn 2$   be
the blowing-up of $\pn 2$ at a point $x$,  with  exceptional curve $e_1$, and set
$L_1:=\sigma_1^*\sO_{\pn 2}(2) - e_1$. Thus $(\scy_1,L_1) =
({\mathbb F}_1, [C_0+2f])$, where $C_0=e_1$ is the minimal section and
$f$ is a fibre of ${\mathbb F}_1$. In particular, $L_1$ is very ample.
The surface $\scy_1$ has two
extremal rays, namely $R'=\reals_+[e_1]$, which is not nef, and $R''=\reals_+[f]$, which is nef, and both are 
$L_1$-positive. Note that $\tau_{L_1}(R')=1$ while
$\tau_{L_1}(R'')=2$. We have $tK_{{\mathbb F}_1}+L_1= (1-2t)C_0 +
(2-3t)f$, which is nef if and only if $1-2t \geq 0$ and $2-3t\geq
1-2t$. Therefore $\frak{t}(\scy_1, L_1)= \frac{1}{2}$. Hence the
ray orthogonal to $\frak{t}(\scy_1,L_1)K_{\scy_1}+L_1$ is
$R''$. 

Now let $\sigma: \sM \to \scy_1$ be the blowing-up of ${\mathbb F}_1$ at a point
$x_1$ lying on $e_1$, let  $e$ be the exceptional curve, and set $\sL=
\sigma^*L_1-e$, which is nef and big, in fact spanned (see e.g., \cite[Lemma 1.7.7]{Book}).
Note that $\sM$
contains exactly two $(-1)$-curves, namely $e$ and
$\widetilde{f_0}$, the proper transform of the fiber $f_0$ of
${\mathbb F}_1$ containing $x_1$. In fact, $R:=\reals_+[e]$ and
$R_0:=\reals_+[\widetilde{f_0}]$ are the only two extremal rays. 
We have $\sL \cdot e =
1$, while
$\sL \cdot \widetilde{f_0}=0.$ Hence $R$ is $\sL$-positive, with 
$\tau_{\sL}(R) = 1$, while $R_0$ is not,
so that
$(\sM,\sL)$ is not rays-positive.  

The
 first reduction of $(\sM, \sL)$ is $(M,L)=(\pn 2, \sO_{\pn 2}(2))$ with first reduction map $\Phi= \sigma_1 \circ
\sigma$, whose exceptional locus is $e \cup \widetilde{e_1}$. In
conclusion, 
the first reduction $(M,L)$ is as in case $(4)$ of Proposition \ref{highnef}. By the
way note that $(M,L)$ is also the first reduction of $(\scy_1,
L_1)$, which is not defined in the classical case.
\end{example*}

Here are some   examples of rays-positive manifolds $(\sM,\sL)$. The first two of them 
fit  in case $(3)$ of Proposition \ref{highnef}. 
\begin{example*}\label{really} Let  $\scv_n$ be a degree $1$ indecomposable vector bundle of rank $n$ over a smooth curve $Y$ of genus  $1$. It is well-known that $\scv_n$ is  ample for any $n\geq 1$. Hence the tautological line bundle  of $\scv_n$ is ample on ${\mathbb  P}(\scv_n)$.

Now, let $\sE:=\sO_Y^{\oplus s} \oplus \scv_{n-s}$, for some positive integer $s$, $\sM:={\mathbb P}(\sE)$, and let  $\sL$ be the tautological bundle  of $\sE$ on $\sM$. Then $\sL$  is nef and $\sL^n = \deg( \sE)= 1$.  Moreover $g(\sM,\sL)=g(Y)=1$.
In particular we get  an example of a quasi-polarized manifold as in  \cite[p. 109]{FujitaQ}.
Moreover, the canonical bundle formula gives
$K_\sM = -n\sL + p^*(K_Y+\scv_1) = -n\sL + F$,  where $p\colon \sM\to Y$ is the bundle projection and   $F$ is a fiber.   Therefore, denoting  by ``$\equiv$'" the numerical equivalence,  one has $K_F \equiv -n\sL_F$. Thus $\sL\cdot C>0$ for any rational extremal curve $C\subset F$, so that $(\sM,\sL)$ is rays-positive according to Lemma  \ref{Dm}.
\end{example*}

\begin{example*}\label{exnup}  Let $C$ be a non-singular curve of genus $g\geq 2$. There exists a stable vector bundle $\sE$ of rank $2$ and degree zero on $C$
 whose tautological line bundle on ${\mathbb P}(\sE)$   is nup and not
big (see \cite[Example 10.6]{HrtAmple}). Let $A$ be an ample line bundle
on $C$ and set $\sM:={\mathbb P}(\sE\oplus A)$. Note that $\sE\oplus A$ is not ample since  $\sE$ has degree zero. Then the tautological line bundle $\sL$ on $\sM$ is not ample, but nup and big \cite[3.13]{BSnup}, whence rays-positive.
\end{example*}

\begin{example*}\label{quasiF}Let $\sM$ be  a smooth projective variety whose anticanonical  bundle $-K_\sM$  is nef and not numerically trivial. Clearly $\ft(\sM,-K_\sM)= 1$, hence $(\sM,-K_\sM)$ is  rays-positive. In particular, if $\sM$ is a quasi-Fano manifold of index $r$ and $-K_\sM=r\sL$ for a nef and big line bundle $\sL$, then $\ft(\sM,\sL)=\frac{1}{r}$,  so that  $(\sM,\sL)$ is rays-positive. Note that if $\sL$ is nup then it is ample  by the basepoint free theorem (cf.  Errata to \cite{Debarre}, p.  219). This shows that $\sL$ is rays-positive but not nup whenever  $(\sM,\sL)$ is quasi-Fano but not Fano.

A general construction is as follows. Let $(Y,H)$ be an $m$-dimensional quasi-Fano manifold
of index $r$, $-K_Y=rH$. Set $\sM:={\mathbb P}(\sE)$, where $\sE= \sO_Y^{\oplus s} \oplus H^{\oplus r}$
and let $\xi$ be the tautological line bundle on $\sM$.
Then the canonical bundle formula yields $K_{\sM}=-(r+s)\xi$. Note that $\xi$ is nef, so being $\sE$. If it is also big,
then $(\sM, \xi)$ is a quasi-Fano manifold of index $r+s = \dim \sM + 1 - m$.
For instance, take $m=3$. Then $1\leq r\leq 4$ and $s+3 \leq \dim \sM \leq s+6$ accordingly. 
The Chern polynomial of $\sE$ is $c(\sE;t)= (1+Ht)^r\ \text{mod}\ H^4$.  Recall that $H^3 >0$. 
By an iterated application of the Chern--Wu formula we have 
$\xi^{\dim \sM} = \xi^{r+s+2} =  c_1^3 -2c_1c_2 + c_3 > 0$, where $c_i=c_i(\sE)$,
showing that $\xi$ is big. Therefore, regardless the value of $r$,  $\sM$ is a quasi-Fano manifold, of index $r+s= \dim \sM - 2$.

Examples of quasi-Fano threefolds of this type are discussed in \cite[Example 2.10 and Proposition 3.2]{CJR}.
\end{example*}

Further examples of rays-positive  manifolds come from projective varieties  with mild singularities. The following definition and result will find an application in \S \ref{lowdeg}. 

\begin{definition*} \label{def:crepsing}
 Let $X$ be a reduced and irreducible variety. We say that $X$ has  {\em crepant singularities} if the normalization $X'$ of $X$ is $\mathbb{Q}$-Gorenstein, that is,  the canonical Weil divisor $K_{X'}$ is $\mathbb{Q}$-Cartier, and $X'$ admits a resolution of singularities $\rho:\widetilde{X} \rightarrow X'$ such that $K_{\widetilde{X}}=\rho^*K_{X'}$. (Clearly, smooth varieties have crepant singularities.)
We say that the composition morphism $\pi:\widetilde{X} \rightarrow X$ is a {\em crepant resolution} of $X$.
\end{definition*}

\begin{lemma} \label{crep-am}
 Let $X$ be a variety with crepant singularities and $L$ an ample line bundle on $X$. Let $\pi:\widetilde{X} \rightarrow X$ be any crepant resolution.
Then $(\widetilde{X}, \pi^*L)$ is rays-positive.
\end{lemma}
\proof
If $(\widetilde{X},\pi^*L)$ is not rays-positive, then  by Lemma \ref{Dm} there is a curve $C$ such that $K_{\widetilde{X}}\cdot C<0$ and $\pi^*L\cdot C=0$. 
Since $L$ is ample, this  means that $C$ is contracted by $\pi$. Let
$\pi:\widetilde{X} \stackrel{\rho}{\longrightarrow} X' \stackrel{\nu}{\longrightarrow} X$
be the Remmert--Stein factorization of $\pi$, where
$\nu:X' \rightarrow X$ is the normalization. Then $C$ is contracted by $\rho$ and, since  $\pi:\widetilde{X} \rightarrow X$ is a crepant resolution, we have
$K_{\widetilde{X}}\cdot C =  \rho^*K_{X'} \cdot C = K_{X'}\cdot \rho(C)=0$, 
a contradiction.
\qed

 \section{Structure results  for rays-positive manifolds}\label{FRPM}
 \addtocounter{subsection}{1}\setcounter{theorem}{0}
 
In this section we get some classification results for   rays-positive manifolds $(\sM,\sL)$.  Note that Proposition \ref{highnef} already classifies such pairs with $\ft(\sM,\sL) < \frac{1}{n-1}$ (equivalently, with  $K_{\sM}+(n-1)\sL$ non-nef)  since  in this case
   there exists an extremal ray $R$ such that $\tau_\sL(R)>n-1$. In view of Remark \ref{remnefvalue}, this is in fact  a  consequence of what is known for  polarized manifolds.

 For higher values of $\ft(\sM,\sL)$ the first reduction enters in the picture.
 We start with the following result (cf. the classical case where $\sL$ is ample).

\begin{lemma} \label{prop:cara} Let $(\sM,\sL)$ be a rays-positive manifold
 of dimension $n\geq 2$ such that $K_{\sM}$ is not nef. Let $(M,L)$ be a first reduction of $(\sM,\sL)$.
Then the following conditions are equivalent:
\begin{enumerate}
\item[$(1)$] $K_{M}$ is not nef and $\ft(M,L) \leq \frac{1}{n-1}$.
\item[$(2)$] $(M,L)$ is as in one of  propositions $\ref{highnef}$ and  $\ref{notbig}(1)$--$(3)$.
\end{enumerate}
Moreover, $(1)$ and  $(2)$  imply
\begin{enumerate} 
\item[$(3)$] $K_{\sM}+(n-1)\sL$ is not  nef and big.
\end{enumerate}
If, furthermore, $\sL$ is big, then condition $(3)$ is equivalent to  $(1)$ and $(2)$. 
\end{lemma}
\proof By Lemma  \ref{tsale} the pair $(M,L)$  is rays-positive. A direct check shows that  $(2)$ implies $(1)$. The converse follows since case $(4)$ of 
Proposition \ref{notbig} cannot occur on $M$ by definition of  first reduction. Therefore, $(1)$ and $(2)$ 
are equivalent and one easily sees that $K_{M}+(n-1)L$ is not nef and big in these cases, whence nor is
$K_{\sM}+(n-1)\sL$ as $K_\sM+(n-1)\sL=\Phi^*(K_M+(n-1)L)$. 

In view of Remark \ref{remnefvalue}, the fact that $(3)$ implies $(1)$ if 
$\sL$ is big follows from classical adjunction.
\qed

Now if one is interested in a biregular, and not only birational, classification of varieties, an interesting question to ask is whether the first reduction map $\Phi:\sM \rightarrow M$ is an isomorphism or not in the equivalent conditions $(1)$ and  $(2)$ of Lemma \ref{prop:cara} (recall that  by Lemma \ref{tsale}, the map $\Phi$ is an isomorphism if $\ft(\sM,\sL)>\frac{1}{n-1}$).
The next two results, which will be proved together,  deal with this question in the cases where $\sL$ is, respectively, rays-positive and ample.  The ample case, treated in Proposition \ref{nonid-am},  is included to make the comparison  with the classical case of polarized manifolds.

\begin{prop} \label{nonid-rp}  Let $(\sM,\sL)$ be a rays-positive manifold
 of dimension $n\geq 2$ and $(M,L)$  a first reduction of $(\sM,\sL)$ with first reduction map 
$\Phi:\sM \rightarrow M$. Assume that $K_M$ is not nef and $\ft(M,L) \leq \frac{1}{n-1}$. 
If $\Phi$ is not an isomorphism, then  either
\begin{enumerate}
\item[$(1)$]  $(M,L) = (\mathbb{P}^2, \sO_{\mathbb{P}^2}(2))$, $(\sM,\sL) = (\mathbb{F}_1, [C_0+2f])$ 
and $\Phi$ is the contraction of the $(-1)$-section $C_0$; or
\item[$(2)$] $(M,L)$ is as in Proposition $\ref{notbig}(1)$--$(3)$.
\end{enumerate}
\end{prop}

\begin{prop} \label{nonid-am}
 Let $(\sM,\sL)$ be a polarized manifold of dimension $n\geq 2$ 
and  $(M,L)$  a first reduction of $(\sM,\sL)$ with first reduction map 
$\Phi:\sM \rightarrow M$. Assume that $K_{\sM}+(n-1)\sL$ is not nef and big.
If $\Phi$ is not an isomorphism, then one of the following cases occurs: 
\begin{enumerate}
\item[$(1)$]  $(M,L) =(\mathbb{P}^2, \sO_{\mathbb{P}^2}(2))$, $(\sM,\sL) = (\mathbb{F}_1, [C_0+2f])$ 
and $\Phi$ is the contraction of the $(-1)$-section $C_0$. 
\item[$(2)$]  $n=2$, $\sM$ is a non-minimal Del Pezzo surface and $\sL=-K_{\sM}$.  Equivalently, there is a birational morphism $\theta: \sM \rightarrow \mathbb{F}_1$ expressing $\sM$ as $\mathbb{F}_1$ blown up at $s$ points lying on distinct fibers, $0 \leq s \leq 7$, and $(\theta_*\sL)^{**}=-K_{\mathbb{F}_1}$. Here $(M,L)=(\mathbb{P}^2, \sO_{\mathbb{P}^2}(3))$ and $\Phi$ equals the composition of $\theta$ with the contraction of the section $C_0$. 
\item[$(3)$]  $n=2$, $(\sM,\sL)$ is a conic fibration over a smooth curve $Y$ admitting some reducible fibers and $(M,L)$ is a conic fibration over $Y$ with  irreducible fibers. Here $\Phi$ is the contraction of one component of each reducible fiber.  
\item[$(4)$] $n=3$ and $(\sM,\sL)$ is the Del Pezzo threefold of degree $7$. Equivalently, $\sM=\mathbb{P}(\sO_{\mathbb{P}^2}(2) \oplus \sO_{\mathbb{P}^2}(1))$ and $\sL$ is the tautological line bundle. In this case $(M,L)=(\mathbb{P}^3, 
\sO_{\mathbb{P}^3}(2))$ and $\Phi$ is the contraction of the $(-1)$-plane $E \subset \sM$ representing the tautological section of $\sO_{\mathbb{P}^2} \oplus \sO_{\mathbb{P}^2}(-1)$. 
\end{enumerate}
\end{prop}

\noindent {\it Proofs of Propositions $\ref{nonid-rp}$  and $\ref{nonid-am}$.} 
We work with the assumptions as in Proposition \ref{nonid-rp} (weaker than those in Proposition \ref{nonid-am}).  If $\Phi$ is not an isomorphism, it  factors as
$
\Phi:  \sM \stackrel{\theta}{\longrightarrow}M_x \stackrel{\grs_x}{\longrightarrow}M, $
where $\grs_x$ is  the blowing-up at a point $x \in M$, and $\theta$ is a sequence of blowing-ups (possibly  an isomorphism). We let $E_x \cong \mathbb{P}^{n-1}$ be the exceptional divisor of
$\grs_x$ and $L_x:= \grs_x^*L-E_x$. 

Pick any extremal ray $R=\mathbb{R}_+[C]$ on $M$  such that $\tau_L(R)=\frac{1}{\ft(M,L)} \geq n-1$. Then $R$ is nef and the pair $(M,L)$ is  described as  in propositions \ref{highnef} and \ref{notbig}$(1)$--$(3)$ by Lemma \ref{prop:cara}.  In particular,  the curves algebraically equivalent to $C$ cover $M$ (and they are all smooth rational curves). Thus we can choose one such  curve $\Gamma$ passing through $x$, and we denote by $\Gamma_x \cong \Gamma$ its strict transform on $M_x$.
Therefore
\begin{equation} \label{eq:azero}
L_x\cdot \Gamma_x =  
\left(\grs_x^*L-E_x\right )\cdot \Gamma_x =   L\cdot \Gamma -E_x \cdot \Gamma_x = L \cdot C- 1 
\end{equation}
and
\begin{equation} \label{eq:ameno}
  K_{M_x} \cdot \Gamma_x=   (\grs_x^*K_M+(n-1)E_x)\cdot\Gamma_x = 
  K_M \cdot \Gamma + (n-1)E_x \cdot \Gamma_x= -\ell(R)+n-1. 
\end{equation}

Now let $\{E_i\}$ be the (possibly empty) set of irreducible exceptional divisors of $\theta$. Then 
$
K_{\sM} = \theta^*K_{M_x} + \sum_i \alpha_i E_i$  and 
$\sL = \theta^*L_{x} - \sum_i \beta_i E_i$,  for some positive integers $\gra_i$, $\grb_i$.
Let $\Delta$ be the strict transform of $\Gamma_x$ on $\sM$. Then, by \eqref{eq:azero} and \eqref{eq:ameno}, we have
\begin{equation} \label{eq:azero1}
\sL\cdot  \Delta  =  (\theta^*L_{x} - \sum_i \beta_i E_i)\cdot \Delta = 
 L_x\cdot \Gamma_x - \sum_i \beta_iE_i\cdot \Delta =  L\cdot C - 1 - \sum_i \beta_i   E_i \cdot\Delta
\end{equation}
and
\begin{equation}\label{eq:ameno1} \begin{split}
 K_{\sM}\cdot \Delta =
K_{M_x}\cdot \Gamma_x + \sum_i \alpha_i   E_i \cdot \Delta = -\ell(R)+n-1 + \sum_i \alpha_i E_i\cdot \Delta.
\end{split} \end{equation}

 Since $(M,L)$ and the contraction of $R$ are as in propositions \ref{highnef} or \ref{notbig}(1)--(3), one can directly check that 
$ L\cdot C=1$ (whence $\tau_L(R)= \ell(R)$), except for the following cases:
\begin{itemize}
\item[(a)] $L\cdot C=2$, $\ell(R)=3$, $(M,L)=(\mathbb{P}^2, \sO_{\mathbb{P}^2}(2))$. 
\item[(b)] $ L\cdot C=\ell(R)=2$, $(M,L)$ as in Proposition  \ref{notbig}$(2)$ with $n=2$. 
\item[(c)] $L\cdot C=\ell(R)=3$, $(M,L)=(\mathbb{P}^2, \sO_{\mathbb{P}^2}(3))$. 
\item[(d)] $ L\cdot C=2$, $\ell(R)=4$, $(M,L)=(\mathbb{P}^3, \sO_{\mathbb{P}^3}(2))$.
\end{itemize}

Before  dealing with  cases (a)--(d), we treat the case when $L\cdot C=1$.

Since $\sL$ is nef, we obtain from \eqref{eq:azero1} that $ \sL\cdot\Delta=0$ and $ E_i\cdot \Delta=0$ for all $i$.  In particular, $\sL$ is not ample, so this case does not occur in Proposition \ref{nonid-am}. 
Inserting into \eqref{eq:ameno1}, we obtain
$ K_{\sM}\cdot \Delta = -\ell(R)+n-1= -\tau_L(R)+n-1 \leq 0$. Since $\sL$ is rays-positive with $\sL\cdot \Delta=0$ we must have $ K_{\sM}\cdot\Delta \geq 0$  by Lemma \ref{Dm}$(3)$.
Therefore
$\tau_L(R) =n-1$, so that $(M,L)$ is as in 
Proposition \ref{notbig}$(1)$--$(3)$. This leads to  case $(2)$ in Proposition \ref{nonid-rp}. 

Now we treat the cases (a)--(d) separately. Only case (a) of them  is not as in Proposition \ref{notbig}$(1)$--$(3)$.  In this case we have $(M_x,L_x)=(\mathbb{F}_1, [C_0+2f])$ and $\grs_x$ is the contraction of the $(-1)$-section $C_0$. Note that $(M_x,L_x)$ is as in 
Proposition \ref{highnef}(3) with extremal ray therein $R'=\reals_+[C']$ satisfying $L_x\cdot C'=1$ and $\tau_{L_x}(R')=n$. 
Repeating the same argument as above with $(M,L)$ replaced by $(M_x,L_x)$ shows that, if $\theta$ were not an isomorphism, then $\tau_{L_x}(R')=n-1$, a contradiction. 
Thus $\theta$ must be an isomorphism, so that we end up  in case $(1)$ of propositions \ref{nonid-rp}  and   \ref{nonid-am}. 
This completes the proof of  Proposition \ref{nonid-rp}. 

To finish the proof of Proposition \ref{nonid-am}, we will now assume that $\sL$ is ample.

In case (b), the pair $(M,L)$ is a conic fibration over a smooth curve $Y$ with irreducible
fibers. Because of the ampleness of $\sL$, it thus follows that   $\Phi$ is  a blowing-up of distinct points on distinct fibers. This yields case $(3)$ of Proposition \ref{nonid-am}. 

In case (c) we must  have $(M_x,L_x)=(\mathbb{F}_1,[2C_0+3f])$. In particular, $L_x=-K_{\mathbb{F}_1}$, and, by the properties of the first reduction map, $\sL=-K_{\sM}$, so that $(\sM, \sL)$ is a Del Pezzo surface.  Again the ampleness of $\sL$ implies that $\Phi$ is a blowing-up of $s$ points, lying on distinct fibers,
and $s <8$ because  $0< \sL^2=K_{\mathbb{F}_1}^2-s=8-s$. This gives case $(2)$ of Proposition \ref{nonid-am}.

In case (d), we have $(M,L)=(\mathbb{P}^3, \sO_{\mathbb{P}^3}(2))$ and 
$(M_x,L_x)$ is $\mathbb{P}(\sO_{\mathbb{P}^2}(2) \oplus \sO_{\mathbb{P}^2}(1))$ with its tautological line bundle, and $\Phi$ is the contraction of the plane representing the tautological section of $\sO_{\mathbb{P}^2} \oplus \sO_{\mathbb{P}^2}(-1)$. To show that we are in case $(4)$ of Proposition \ref{nonid-am}, we must only show that $\theta$ is an isomorphism.  But $(M_x,L_x)$ is as in Proposition \ref{notbig}$(3)$, with $Y=\mathbb{P}^2$, and with extremal ray  therein $R'=\reals_+[C']$  satisfying $L_x\cdot C'=1$.
Therefore, repeating the same argument as above with $(M,L)$ replaced by $(M_x,L_x)$ shows that if $\theta$ were not an isomorphism, then $\sL$ would not be ample, a contradiction.
 (More directly,  since there is there is a line passing 
 through any pair of points in $\mathbb{P}^3$, also infinitely near,  one easily sees that blowing-up $M$ at more than one point  would make $\sL$ not ample.)
\qed

\begin{rem*} \label{casidoppi}  
 Note that in case $(1)$ of propositions \ref{nonid-rp} and \ref{nonid-am}, the pair $(\sM,\sL)$ is already as in Proposition \ref{highnef}$(3)$; and in case $(4)$ of  Proposition \ref{nonid-am}, the pair $(\sM,\sL)$ is already as in  Proposition \ref{notbig}$(3)$ (with $Y=\mathbb{P}^2$).
\end{rem*}

The next  example shows that the result in Proposition \ref{nonid-rp}$(2)$ is optimal, in the sense that all the possibilities of Proposition \ref{notbig}$(1)$--$(3)$ do in fact occur with $\Phi$ not an isomorphism.  It  exhibits    rays-positive line bundles $\sL$ which are not nup, and  first reductions $(M,L)$ covered by lines as well.

\begin{example*}\label{AnCa}  Let $(M,L)$ be as in Proposition \ref{notbig}$(1)$--$(3)$ with $L$ ample and spanned. Let $\grs:\sM\to M$ be the blowing-up at a point $x\in M$, and let $E\cong\pn {n-1}$ be the exceptional divisor. Then $\sL:=\grs^*L-E$   is nef (see e.g.,  \cite[Lemma 1.7.7]{Book}) and $\sL+(n-1)K_\sM=\grs^*(L+(n-1)K_M)$. Thus $\sL+(n-1)K_\sM$ is nef since $L+(n-1)K_M$ is nef being $\ft(M,L)=\frac{1}{n-1}$. Therefore $\ft(\sM,\sL)\geq\frac{1}{n-1}$, so that $\ft(\sM,\sL)=\frac{1}{n-1}$ by Lemma \ref{tsale}. Hence $(\sM,\sL)$ is rays-positive. Clearly, $\grs:\sM\to M$  is the first reduction map.

Furthermore, if we are in the cases where the extremal ray $R=\reals_+[C]$ described in Proposition \ref{notbig}$(1)$--$(3)$  satisfies the condition $L\cdot C=1$ (as noted in the proof of Proposition \ref{nonid-rp} this happens except for the cases (b), (c), (d) listed in that proof), then the strict transform 
$\Gamma$ of any  curve numerically equivalent to $C$ passing through $x$ satisfies the conditions
$K_\sM\cdot \Gamma=\sL\cdot 
\Gamma=0$. In particular, $\sL$ is rays-positive but not nup.

As concrete examples  of pairs $(M,L)$ as above we may take either a Del Pezzo  $n$-fold of degree $d=L^n$, $3\leq d\leq 4$ (in this case $L$ is very ample and $M$ is covered by lines, see \cite[Chapter 1, \S 8]{FuBook}), or $(M,L)=(\sQ\times B, p_1^*\sO_\sQ(1)+p_2^*\sO_B(3b))$, where $\sQ$ is a smooth hyperquadric in $\pn n$ with $n\geq 4$, $B$ is a smooth curve of genus $g(B)=1$, $b$ is a point on $B$, and  $p_1$, $p_2$ are  the projections on the two factors, or 
 $(M,L)=(S\times \pn {n-2},p_1^*A +p_2^*\sO_{\pn {n-2}}(1)) $, where $S$ is a smooth   surface,  $p_1$, $p_2$ are the projections on the two factors, and $A$ is a very  ample line bundle  on $S$, so that $L$ is very ample on $M$ (note that in this case  $(M,L)$ is a scroll over $S$, in both the classical and the adjunction theoretic sense, since $K_M +(n-1)L=(n-1)p_1^*(A)$.)
\end{example*}

Let $(\sM,\sL)$ be a polarized manifold of dimension $n\geq 2$ such that
$K_{\sM}+(n-1)\sL$ is not nef and big. By Lemma \ref{prop:cara} and Proposition \ref{nonid-am} (together with  Remark \ref{casidoppi}), the pair $(\sM,\sL)$ is as in propositions \ref{highnef} and  \ref{notbig}$(1)$--$(3)$, except for the cases $(2)$ and $(3)$ of Proposition \ref{nonid-am}.   These are well-known biregular classification results in adjunction theory (see e.g., \cite[(7.2.1), (7.2.2), (7.2.4), (7.3.2)$(1)$--$(3)$]{Book}).

Note that n the case when  $(\sM,\sL)$  is merely  a rays-positive quasi-polarized manifold the corresponding classification result  is only birational as shown by Example \ref{AnCa}. The precise statement is the following.

\begin{corollary} \label{9ago}
Let $(\sM,\sL)$ be a rays-positive quasi-polarized  manifold of dimension $n\geq 2$. 
If $K_{\sM}+(n-1)\sL$ is not nef and big, then either $(\sM,\sL)$ is as in Proposition $\ref{highnef}$ or  its first reduction is as in Proposition $\ref{notbig}(1)$--$(3)$. \end{corollary}
\proof
Combine Lemma \ref{prop:cara},  Proposition \ref{nonid-rp} and   Remark \ref{casidoppi}. \qed

We also have the following consequence of the results above.

\begin{corollary} \label{cor:nonid-rp} Let $(\sM,\sL)$ be a rays-positive quasi-polarized  manifold of dimension $n\geq 3$. 
If $K_{\sM}+(n-1)\sL$ is not nef and big, then $(\sM,\sL)$ is uniruled of $\sL$-degree at most one 
unless $(\sM,\sL)=(\mathbb{P}^3, \sO_{\mathbb{P}^3}(2))$. \end{corollary}
\proof
First note that if $(M,L)$ is uniruled of $L$-degree at most one, then, as
a consequence of the properties of the first reduction map $\Phi$,
 the same is true for $(\sM,\sL)$. 

By Lemma \ref{prop:cara}, the pair $(M,L)$ is as in propositions \ref{highnef}  and  \ref{notbig}$(1)$--$(3)$. As observed in the proof of Proposition \ref{nonid-rp}, for $n\geq 3$, the only pair  $(M,L)$ not uniruled of $L$-degree at most one is  $(\mathbb{P}^3, \sO_{\mathbb{P}^3}(2))$. In this case, if $\Phi$ is not an isomorphism, then $(\sM,\sL)$ is uniruled of $\sL$-degree at most one.  Indeed, take any point $x \in \mathbb{P}^3$ over which $\Phi$ is not an isomorphism. Then $\mathbb{P}^3$ is covered by the family of lines through $x$. These have degree two with respect to $L$ and their strict transforms have  degree $\leq 1$ with respect to $\sL$. 
\qed

Note that in both corollaries \ref{9ago} and \ref{cor:nonid-rp} the condition that $\sL$ is big (i.e., $(\sM,\sL)$ quasi-polarized) and $K_{\sM}+(n-1)\sL$ is not nef and big can be replaced by the weaker condition that $K_M$ is not nef and $\ft(M,L) \leq \frac{1}{n-1}$, where  $(M,L)$ is the first reduction of $(\sM,\sL)$. This follows from Lemma \ref{prop:cara}.

\section{Pseudo-effectivity and varieties of low degree}
\label{lowdeg}\addtocounter{subsection}{1}\setcounter{theorem}{0}

The aim of this section is to classify projective  varieties with   crepant singularities (see Definition \ref{def:crepsing}) and  of small degree with respect to the codimension.

The following  proposition yields a classification up to first reductions  of  rays-positive manifolds $(\sM,\sL)$ such that $K_{\sM}+(n-2)\sL$  is  not {\em pseudo-effective}, i.e., not contained in  the closure of the cone spanned by classes of effective divisors.

\begin{prop} \label{prop:class} Let $(\sM,\sL)$ be a rays-positive manifold of dimension 
 $n\geq 3$.  If  $K_{\sM}+(n-2)\sL$ is not pseudo-effective,   then either  
$(\sM,\sL)$ is described as in  Proposition  $\ref{highnef}$, or
its  first reduction $(M,L)$  is described as in one  of propositions {\rm \ref{notbig}$(1)$--$(3)$} and  $\ref{n21}$.
\end{prop}
\proof 
By  \cite[Theorem 0.2]{BDPP},   the fact that $K_{\sM}+(n-2)\sL$ is not
pseudo-effective is equivalent to the existence of  a covering family of curves on $\sM$
such that 
$(K_{\sM}+(n-2)\sL) \cdot C <0$ for all curves $C$ in the family. 
In particular $K_\sM$ is not nef. 

Consider the first reduction $(M,L)$  of $(\sM,\sL)$ with first reduction morphism $\Phi$. Then $L$ is nef and rays-positive by Lemma \ref{tsale}.
Recall that $\sL = \Phi^*L -J$, where $J$ is an effective $\Phi$-exceptional divisor. 
The general curve $C$ in the family above
is not contained in the support of $J$, so that
$J \cdot C \geq 0$,
and $C$ is not contracted by $\Phi$. Let $\Gamma:=\Phi(C)$.
Then, since $\Phi^*(K_M+(n-1)L)=K_\sM+(n-1)\sL$, we have
$$
(K_{M}+(n-2)L) \cdot \Gamma  = (K_{\sM}+(n-2)\sL) \cdot C - J \cdot C \leq  (K_{\sM}+(n-2)\sL) \cdot C <0.$$
Hence $K_M$ is not nef and  $\ft(M,L) < \frac{1}{n-2}$.
By Lemma \ref{prop:cara}, either $(M,L)$ is as in one of  propositions \ref{highnef} and \ref{notbig}$(1)$--$(3)$ or $\frac{1}{n-1}<\ft(M,L)<\frac{1}{n-2}$. 
In the latter case,  $(M,L)$ is as in Proposition \ref{n21}. Finally, if $(M,L)$ is as in Proposition \ref{highnef}, then by Proposition \ref{nonid-rp} and Remark \ref{casidoppi},  the pair $(\sM,\sL)$
is  as in Proposition \ref{highnef} as well, concluding the proof. 
\qed

As an application, we extend  the main result in \cite{IoSmall}, providing  a classification of projective  varieties with  crepant singularities  and   small degree. Note that the assumption  $d<2\ {\rm codim}_{\pn N}(X)+2$  in the theorem  below can be rephrased in terms of $\Delta$-genus as $d>2\Delta(X,\sO_X(1))$.

\begin{theorem} \label{thM:lowdeg}
 Let $X \subset \mathbb{P}^N$ be a reduced and irreducible 
variety of dimension $n\geq 3$ and  degree $d$, and with crepant singularities. Assume  $d < 2\ {\rm codim}_{\pn N}(X)+2$. Let $\pi:\sM \rightarrow X$ be any crepant resolution and let $\sL:=\pi^*\sO_X(1)$. Then   either  
$(\sM,\sL)$ is described as in  Proposition  $\ref{highnef}$, or
its  first reduction $(M,L)$  is described as in one  of propositions {\rm \ref{notbig}$(1)$--$(3)$} and  $\ref{n21}$. Moreover, in the scroll case $(3)$ of Proposition $\ref{notbig}$ the base  surface   is  ruled.\end{theorem}
\proof 
We have that $\sL$ is globally generated with $\dim |\sL| \geq \dim |\sO_X(1)|=N$. We can pick $n-1$ general members $H_1, \ldots, H_{n-1}$ in $|\sL|$ such that each $\sM_i : = H_1 \cap \cdots \cap H_i$, with  $i=1, \ldots, n-1$,
is smooth and irreducible of dimension $n-i$. We let $\sM_0=\sM$. From the standard  restriction sequences we get 
$\dim|\sL_{\sM_{i+1}}| \geq \dim|\sL_{\sM_{i}}| -1$,
so that $\dim|\sL_{\sM_{i}}| \geq N-i$. In particular, on the smooth curve $C := \sM_{n-1}$ we have, by assumption, 
\[ \deg (\sL_C)-2 \dim |\sL_C| \leq d -2(N-(n-1))= d-2\ {\rm codim}_{\pn N}(X)-2<0. \]
 Thus by Clifford's theorem we must have $h^1(\sL_C)=0$, so that
\[ \chi(\sL_C) =h^0(\sL_C) \geq N -(n-1)+1 =N-n+2. \]
Consider the smooth surface $S:=\sM_{n-2}$. By  the Riemann--Roch theorem we get 
\[   K_S\cdot C=
C^2 -2\big(\chi(\sL_S) - \chi(\sO_S)\big)
=C^2-2\chi(\sL_C) \leq d-2(N-n+2) < -2. \]
Therefore $(K_{\sM}+(n-2)\sL)\cdot C=K_S\cdot C <-2$, whence
$K_{\sM}+(n-2)\sL$ is not
pseudo-effective by \cite[Theorem 0.2]{BDPP}. 
Then  the result follows from Proposition
\ref{prop:class}, $(\sM,\sL)$ being rays-positive by Lemma \ref{crep-am}.
As to the last assertion, note that in case $(3)$ of Proposition \ref{notbig} the scroll projection  maps the ruled surface $S$ surjectively onto the base $Y$.
\qed

Notation as in Theorem \ref{thM:lowdeg}. If we assume $X$ to be smooth, so that $\pi$ is the identity map and $\sL$ is very ample, then by Proposition \ref{nonid-am}  and Remark \ref{casidoppi} the conclusion would be that either $(\sM,\sL)=(X,\sO_X(1))$  is  as in one of propositions \ref{highnef}
and \ref{notbig}$(1)$--$(3)$ or its first reduction is as in Proposition \ref{n21}. This is just Ionescu's result. The new occurrences  in the case of  crepant singularities are therefore precisely the cases where $(M,L)$ is as in Proposition \ref{notbig}$(1)$--$(3)$ with the reduction map $\Phi$ not the identity. The following example   shows that these cases indeed occur with  $X$ singular  having crepant singularities.

 \begin{example*}\label{Crecon} Let $(M,L)$  and $R=\reals_+[C]$   be as in Proposition \ref{notbig}$(1)$--$(3)$ with the additional assumption that   $L$  be  very ample and  $L\cdot C=1$. Let $\grs:\sM\to M$ be the blowing-up at a point $x\in M$, and denote by  $E\cong\pn {n-1}$  the exceptional divisor. Then $\sL:=\grs^*L-E$   is spanned and big, and, as proved in Example \ref{AnCa}, rays-positive. Let $\pi:\sM\to \pn N$ be the generically finite morphism defined by $|\sL|$, and set $X:=\pi(\sM)$. Note that $X$ is the variety obtained by projecting $M$, embedded  by $|L|$, from the point $x$.
 
We claim that {\em  $X$ is not smooth and $\pi:\sM\to X$ is a crepant resolution of $X$.}
 
 Let us  first show that an (irreducible) curve $\Gamma$ is contracted by $\pi$ if and only if $\Gamma$ is  the  strict transform  under $\grs$ of a line  on $M$ passing through $x$.

Indeed,  let $\ell$ be such a strict transform. Then $\sL\cdot \ell=\grs^* L\cdot \ell-E\cdot \ell=1-1=0$, whence $\pi(\ell)$ is a point.  As to the converse, note that a curve    $\Gamma$  contracted by $\pi$ is not contained in $E$, 
since  $\sL_E\cong\sO_{\pn {n-1}}(1)$. Then $\gamma:=\grs(\Gamma)$ is an   irreducible curve  in $M$, and one has
$0=\sL\cdot\Gamma=\grs^*L\cdot \Gamma-E\cdot\Gamma=L\cdot \gamma-{\rm mult}_x(\gamma)$.  Therefore $\gamma$ (embedded by $|L|$) is an irreducible curve  of degree $d$ with a singular point of multiplicity $d$.  Thus  $d=1$, that is, $\gamma$ is a line passing through $x$, showing the desired assertion.

Since all the curves $\Gamma$ contracted by $\pi$   are  
 strict transforms under $\grs$ of lines on $M$ passing through $x$, they satisfy $\sL\cdot \Gamma=K_\sM\cdot \Gamma=0$ and $E\cdot \Gamma=1$. It is then a standard fact that $X$ is singular (see e.g., \cite[Proposition 1.45]{Debarre}).
 
  Let $\Delta\subset\sM$ be the locus covered by  the  curves $\Gamma$. Since  $E\cdot \Gamma=1$,  any such curve intersects $E$ in precisely  one point. As $\pi_{|E}$ is an isomorphism, we then infer  that
\begin{equation}\label{sestri1}
\pi(\Delta)=\pi(\Delta\cap E)\cong \Delta\cap E,
\end{equation} 
and  that $\pi_{|\Delta'}:\Delta' \to \pi(\Delta'\cap E)$ is a $\pn 1$-fibration for every irreducible component $\Delta'$ of $\Delta$. One has ${\rm codim}_\sM(\Delta)\geq 2$ by \cite[(11.13)]{FuBook}.

  Consider the Remmert--Stein factorization 
$\sM \stackrel{\pi_1}{\longrightarrow} X'\stackrel{\pi_2}{\longrightarrow} X$ of $\pi$. We have  proved that $\Delta={\rm Exc}(\pi_1)$.
 
As shown in Example \ref{AnCa},  we have $\ft(\sM,\sL)=\frac{1}{n-1}$, so that $K_\sM+(n-1)\sL$ is nef. Therefore $m(K_\sM+(n-1)\sL)$ is spanned for $m\gg 0$. Let $f$ be the morphism defined by $|m(K_\sM+(n-1)\sL)+\sL |$ and set $Z:=f(\sM)$. Consider the Remmert--Stein factorization $\sM \stackrel{f_1}{\longrightarrow} X''  \stackrel{f_2}{\longrightarrow} Z$ of $f$.  We may assume that there exists an integer $k\gg 0$ such that the complete linear systems $|k\sL|$ and $|k(m(K_\sM+(n-1)\sL) + \sL)|= 
|k\sL+ km(K_\sM + (n-1)\sL)|$ define the morphisms $\pi_1$  and $f_1$, respectively. Since
$m(K_\sM+(n-1)\sL)$ is spanned, we have a factorization $\pi_1:\sM\stackrel{f_1}{\longrightarrow}X'' \stackrel{j}{\longrightarrow}X'$.

If $\Gamma$ is a curve contracted by $\pi_1$, then $\sL\cdot \Gamma=K_\sM\cdot \Gamma=0$, so that $\Gamma$ is also contracted by $f_1$. 
It thus follows that $j$ is an isomorphism.
Therefore there are Cartier divisors $D_1$ and $D_2$ on $X'$ such that  $\sL=\pi_1^*( D_1)$ as well as $m(K_\sM+(n-1)\sL)+\sL=\pi_1^*(D_2)$. Hence 
$$
mK_\sM = m(K_\sM+(n-1)\sL)+\sL-(m(n-1)+1)\sL
= \pi_1^*\big(D_2-(m(n-1)+1)D_1\big).
$$
Thus $K_\sM=\pi_1^*(\scd)$  for some $\rat$-Cartier divisor $\scd$ on $X'$.

Now let $\omega_{X'}$ be the canonical sheaf on $X'$. It is a reflexive rank $1$ sheaf  defined by $\iota_*\omega_{{\rm Reg}(X')}$, where $\iota:{\rm Reg}(X')\hookrightarrow X'$ is the inclusion of the smooth points.  
Denote by $K_{X'}$ the corresponding Weil divisor (cf. e.g., \cite[Proposition 5.75]{KM}). On the Zariski open set $\sM\setminus{\rm Exc}(\pi_1)$  the strict transform ${\pi_1}_*^{-1}(K_{X'})$  and $K_\sM$ agree. Hence  they agree on $\sM$, as ${\rm codim}_\sM({\rm Exc}(\pi_1))\geq 2$. Thus 
${\pi_1}_*^{-1}(K_{X'}) = K_\sM=\pi_1^*(\scd)$. By pushing down cycles under $\pi_1$, we obtain
$K_{X'}=\scd$. Therefore $K_{X'}$ is a $\rat$-divisor on $X'$ and $K_\sM={\pi_1}^*(K_{X'})$.  This proves  that $\pi$ is a crepant resolution.

\end{example*}

\section{On the sectional genus of rays-positive manifolds}
\label{FC}\addtocounter{subsection}{1}\setcounter{theorem}{0}

As a final application  we prove a special case of a  conjecture of Fujita, and we describe rays-positive manifolds with sectional genus zero or one.

  For  a nef and big line bundle $\sL$ on an $n$-dimensional manifold $\sM$, Fujita \cite{FujitaQ}
conjectured that  $g(\sM,\sL)\geq 0$ and proved  it for $n\leq 3$ by using Mori's results (see \cite[Corollary 4.8]{FujitaQ}). 
We prove the conjecture for   a  rays-positive manifold $(\sM,\sL)$.

Let us note  first that the mere nefness of $\sL$ is  not enough to grant that $g(\sM,\sL)\geq 0$. For instance, let $(\sM,\sL)=(\pn {n-1}\times\pn 1,\sO(b,0))$ with $b\geq 2$. Clearly, $\sL$ is nef and not big. Moreover, the pair $(\sM,\sL)$ is not rays-positive: actually, for every curve $C$ lying in the first factor we have $K_\sM\cdot C<0$ and $\sL\cdot C=0$. The genus formula gives $g(\sM,\sL)=1-b^{n-1}<0$.

\begin{prop}\label{genus} Let $(\sM,\sL)$ be a rays-positive manifold of dimension 
 $n\geq 2$.  Then $g(\sM,\sL)\geq 0$, with equality if and only if $(\sM,\sL)$ is one of the following:
$(\pn n,\sO_{\pn n}(1))$,  $(\sQ,\sO_\sQ(1))$, with $\sQ$ a hyperquadric in $\pn {n+1}$,  $(\pn 2,\sO_{\pn 2}(2))$, or a scroll over a smooth rational curve.
 \end{prop}
\proof Consider first the case when  $K_\sM+(n-1)\sL$ is nef. Since nef line bundles are limit of ample line bundles, it follows that
$2g(\sM,\sL)-2=(K_\sM+(n-1)\sL)\cdot \sL^{n-1}\geq 0$, whence  $g(\sM,\sL)\geq 1$.
Thus we can assume that
$K_\sM+(n-1)\sL$ is not nef. Then $\ft(\sM,\sL) <\frac{1}{n-1}$,
so that   $(\sM,\sL)$ is described as in Proposition \ref{highnef} as mentioned at the beginning of Section \ref{FRPM}.  A direct check shows that the sectional genus is zero  in cases
$(1)$, $(2)$  and $(4)$ of that proposition, while, in  case
 $(3)$,  $g(\sM,\sL)=g(Y)\geq 0$.
\qed

 In dimension $n=3$,  there is a complete classification of  quasi-polarized varieties $(\sM,\sL)$ with sectional genus $g(\sM,\sL)=0,1$ (see \cite[\S 4]{FujitaQ}). 
Recall also that in \cite{BSg2} a complete classification of quasi-polarized $a$-minimal (rays-positive in our terminology)  Gorenstein surfaces $(\sM,\sL)$ with  $g(\sM,\sL)=2$ is worked out. Moreover, in \cite{Maeda}, surfaces $(\sM,\sL)$ with $\sL$ merely nef are classified  for $g(\sM,\sL)=0,1$.

\smallskip

As to the case of sectional genus  $g(\sM,\sL)=1$ we have the following (cf. \cite[(5.4)]{FujitaQ}).

\begin{prop}\label{g01} Let $(\sM,\sL)$ be a rays-positive manifold of dimension  $n\geq 2$. 
Assume that $g(\sM,\sL)=1$. Then   either 
$\sL^n=K_\sM\cdot\sL^{n-1}=0$, or $(\sM,\sL)$ is a quasi-Del Pezzo manifold, or $(\sM,\sL)$ is   a scroll over a smooth elliptic  curve, with $\sL$ big.
\end{prop} 
\proof If either  $\sL$ is not big or  $K_\sM+(n-1)\sL$  is trivial we find the first two cases.
 If $K_\sM+(n-1)\sL$ is not nef but $\sL$ is big, then the same argument as in the proof of Proposition  \ref{genus} and a direct check lead to the third case.
 
In the remaining cases $\sL$ is big and  
$K_\sM+(n-1)\sL$ is nef and non-trivial, so $m(K_{\sM}+(n-1)\sL)$ is spanned and non-trivial for $ m \gg 0$ by  the Kawamata--Shokurov basepoint free theorem. Thus $m(K_{\sM}+(n-1)\sL) \cdot \sL^{n-1} >0$, as $\sL$ is big, whence $g(\sM,\sL)>1$.\qed

The following example  shows 
  that the first case in the proposition above   really occurs. Moreover, it also shows that the inequality $g(\sM,\sL)\geq h^1(\sO_\sM)$, conjectured in the setting of quasi-polarized varieties, is not true dropping the bigness assumption.

\begin{example*}\label{Ascroll} Let $Y$ be a smooth abelian surface, $\sM={\mathbb P}({\sE})$, where $\sE=\sO_Y\oplus\sO_Y$, and let $\sL$ be the tautological line bundle of $\sE$ on $\sM$. Clearly $\sL$ is nef, since $\sE$ is trivial, and $\sL^3=K_\sM\cdot \sL^2=0$. The only extremal ray of $\sM$ is $R=\reals_+[f]$, where $f$ is a fiber of the bundle projection. One has $K_\sM\cdot f=-2$, $\sL\cdot f=1$, so that $\sL$ is rays-positive. \end{example*}

\smallskip

{\noindent{\em Acknowledgments.} The second author would like to thank INdAM for partial support, while the remaining authors
acknowledge partial support of MIUR provided in the framework of PRIN ``Algebraic Geometry etc." (Cofin 2006 and 2008) and  the University of Milano (FIRST 2006-2007) for making this collaboration possible. The authors are grateful to the referee for useful remarks.

\bigskip

\noindent M.C. Beltrametti,
Dipartimento di Matematica, 
Universit\`a di Genova, 
Via Dodecaneso 35, 
I-16146 Genova, Italy. e-mail {\tt beltrame@dima.unige.it}

\smallskip

 \noindent A.L. Knutsen
Department of Mathematics,
University of Bergen,
Johannes Bruns gate 12,
N-5008 Bergen, Norway.
e-mail {\tt andreas.knutsen@math.uib.no}

\smallskip

\noindent
 A. Lanteri,
 Dipartimento di Matematica ``F. Enriques'',
 Universit\`a degli Studi di Milano,
 Via C. Saldini 50,
 I-20133 Milano, Italy.
 e-mail {\tt antonio.lanteri@unimi.it}

\smallskip

\noindent
 C. Novelli, Dipartimento di Matematica Pura ed Applicata,
Universit\`a degli Studi di Padova,
Via Trieste 63,
I-35121 Padova, Italy.
e-mail {\tt novelli@math.unipd.it}

\end{document}